\documentclass[11pt,twoside,letterpaper]{article} %% The same for {book}

\usepackage{pslatex}
\usepackage{fancyhdr}
\usepackage[dvips]{graphicx}
\usepackage{geometry}

\usepackage{amsmath,amsfonts,amssymb,amsthm}
\usepackage{cite}
\usepackage{times,fancyhdr}
\usepackage{vatola}

\theoremstyle{plain}

\theoremstyle{definition}

\theoremstyle{remark}

\def\q{\quad}
\def\qq{\qquad}
\def\qtq#1{\q\t{#1}\q}
 \TagsOnRight
\def\({\left(}
\def\){\right)}
\def\[{\left[}
\def\]{\right]}
\def\t{\hbox}
\def\mod#1{\ (\text{\rm mod}\ #1)}
\def\f{\frac}
\def\e{\equiv}
\def\b{\binom}
\let \pro=\proclaim
\let \endpro=\endproclaim
\begin{document}

\title{
{\begin{flushleft}
%\vskip 0.45in
%{\normalsize\bfseries\textit{Chapter~6}}
\end{flushleft}
%\vskip 0.45in
\bfseries\scshape Tur\'an's Problem for Trees} }
\author{\bfseries\itshape Zhi-Hong Sun$^1$
\thanks{E-mail address: zhihongsun@yahoo.com;
 Website: {\tt http://www.hytc.edu.cn/xsjl/szh}}
 \  and Lin-Lin Wang$^2$
 \thanks{E-mail address:  wanglinlin$
 \t{\underline{\q }}$1986@yahoo.cn}\\
$^1$\,School of Mathematical Sciences, Huaiyin Normal
University\\
Huaian, Jiangsu 223001, People's Republic of
China\\
$^2$\,Center for Combinatorics, Nankai University\\
Tianjin 300071, People's Republic of China
}
%%%%%%%%%%

\date{}
\maketitle
\thispagestyle{empty}
\setcounter{page}{1}

% ------- [First Page Running Head] - place it immediately after title! ------
\thispagestyle{fancy}  \fancyhead[L]{J.  Comb. Number Theory
3(2011), no.1, 51-69} \fancyhead[R]{ } \fancyfoot{}
\renewcommand{\headrulewidth}{0pt}
%------------------------------------------------------------------------------

\begin{abstract}
 For a forbidden graph
$L$, let $ex(p;L)$ denote the maximal number of edges in a simple
graph of order $p$ not containing $L$. Let $T_n$ denote the unique
tree on $n$ vertices with maximal degree $n-2$, and let
$T_n^*=(V,E)$ be the tree on $n$ vertices with
 $V=\{v_0,v_1,\ldots,v_{n-1}\}$ and
 $E=\{v_0v_1,\ldots,v_0v_{n-3},v_{n-3}v_{n-2},v_{n-2}v_{n-1}\}$.
In the paper we give exact values of $ex(p;T_n)$ and $ex(p;T_n^*)$.

%\vspace{.08in}
%\noindent \textbf{Key words and phrases:}

\vspace{.08in}
\end{abstract}

%[Float at the top formula]
\makeatletter \setlength\@fptop{0\p@} \makeatother

% [Clears Header Style on the Last Empty Odd pages]
\makeatletter
\def\cleardoublepage{\clearpage\if@twoside \ifodd\c@page\else%
    \hbox{}%
    \thispagestyle{empty}%
    \newpage%
    \if@twocolumn\hbox{}\newpage\fi\fi\fi}
\makeatother

\renewcommand{\thesection}{\arabic{section}.}
\renewcommand{\thesubsection}{\thesection\arabic{subsection}.}
\renewcommand{\thesubsubsection}{\thesubsection\arabic{subsubsection}.}

\def\figurename{Figure}
\makeatletter
\renewcommand{\fnum@figure}[1]{\figurename~\thefigure.}
\makeatother

\def\tablename{Table}
\makeatletter
\renewcommand{\fnum@table}[1]{\tablename~\thetable.}
\makeatother

\section{Introduction}
In the paper, all graphs are simple graphs. For a graph $G=(V(G),E(G))$
let $e(G)=|E(G)|$ be the number of edges in $G$ and let
$\Delta(G)$ be the maximal degree of $G$.
  For a family of forbidden graphs $L$, let $ex(p;L)$
denote the maximal number of edges in a graph of order $p$ not
containing any graphs in $L$. The corresponding Tur\'an's problem is
to evaluate $ex(p;L)$. For a graph $G$ of order $p$, if $G$ does not
contain any graphs in $L$ and $e(G)=ex(p;L)$, we say that $G$ is an
extremal graph. In the paper we also use $Ex(p;L)$ to denote the set
of extremal graphs of order $p$ not containing any graphs in $L$.
\par Let $\Bbb N$ be the set of
positive integers. Let $p,n\in\Bbb N$ with $p\ge n\ge 2$. For a
given tree $T$ on $n$ vertices,
 it is difficult to determine the value of $ex(p;T)$.
The famous Erd\"os-S\'os
 conjecture asserts that
$ex(p;T)\le \f{(n-2)p}2$. For the progress on the Erd\"os-S\'os
 conjecture, see [2,6,7,8].
 Write
$p=k(n-1)+r$, where $k\in\Bbb N$ and $r\in\{0,1,\ldots,n-2\}$. Let
$P_n$ be the path on $n$ vertices. In [3] Faudree and Schelp showed
that
$$ex(p;P_n)=k\binom {n-1}2+\binom r2.\tag 1.1$$
In the special case  $r=0$, (1.1) is due to Erd$\ddot{\t{\rm o}}$s
and Gallai [1]. Let $K_{1,n-1}$ denote the unique tree on $n$
vertices with $\Delta(K_{1,n-1})=n-1$, and let
  $T_n$ denote the unique tree on $n$ vertices
with $\Delta(T_n)=n-2$. In Section 2 we determine $ex(p;K_{1,n-1})$,
and in Section 3 we obtain the exact value of $ex(p;T_n)$.
\par For $n\ge 4$ let $T_n^*=(V,E)$ be the tree on $n$ vertices with
 $V=\{v_0,v_1,\ldots,v_{n-1}\}$ and
 $E=\{v_0v_1,\ldots,v_0v_{n-3},v_{n-3}v_{n-2},v_{n-2}v_{n-1}\}$.
In Section 4 we completely determine the value of $ex(p;T_n^*)$.

In addition to the above notation, throughout the paper we
also use the following notation: $[x]\f{\q}{\q}$the greatest integer
not exceeding $x$, $d(v)\f{\q}{\q}$the degree of the vertex $v$ in a
graph, $\Gamma(v)\f{\q}{\q}$the set of vertices adjacent to the
vertex $v$, $d(u,v)\f{\q}{\q}$the distance between the two vertices
$u$ and $v$ in a graph, $K_n\f{\q}{\q}$the complete graph on $n$
vertices, $K_{m,n}\f{\q}{\q}$the complete bipartite graph with $m$
and $n$ vertices in the bipartition, $G[V_0]\f{\q}{\q}$the subgraph
of $G$ induced by vertices in the set $V_0$, $G-V_0\f{\q}{\q}$the
subgraph of $G$ obtained by deleting vertices in $V_0$ and all edges
incident with them, $G-M\f{\q}{\q}$the graph obtained by deleting
all edges in $M$ from the graph $G$, $G+M\f{\q}{\q}$the graph
obtained by adding all edges in $M$ from the graph $G$.

 \section{The Evaluation of $ex(p;K_{1,n-1})$}
%   Theorem 2.1
%\begin{theo}
\pro{Theorem 2.1}
Let $p,n\in\Bbb N$ with $p\geq n-1\geq 1$. Then
$ex(p;K_{1,n-1})=[\f{(n-2)p}2]$.
\endpro
%   Proof.
\begin{proof}
Clearly
$ex(n-1;K_{1,n-1})=e(K_{n-1})=\f{(n-1)(n-2)}2$. Thus the result is
true for $p=n-1$. Now we assume $p\ge n$.
 Suppose that $G$ is
a graph of order $p$ without $K_{1,n-1}$. Then clearly
$\Delta(G)\leq n-2$ and so $2e(G)=\sum_{v\in V(G)}d(v)\leq
p\Delta(G)\leq (n-2)p$. Hence, $ex(p;K_{1,n-1})\leq \f{(n-2)p}2$.
As $ex(p;K_{1,n-1})$ is an integer, we have
$$ex(p;K_{1,n-1})\leq\big[\f{(n-2)p}2\big].\tag 2.1$$
%\end{theo}

Clearly $ex(p;K_{1,1})=0$. So the result holds for $n=2$.
As $[\f p2]K_2$ does not contain $K_{1,2}$, we have
$ex(p;K_{1,2})\geq [\f p2]$. This together with (2.1) gives
$ex(p;K_{1,2})=\[\f p2\]$. So the result is true for $n=3$.
\par Suppose that $G$ is a Hamilton cycle with $p$
vertices. Then $G$ does not contain $K_{1,3}$. Thus we have
$ex(p;K_{1,3})\geq p$. Combining this with (2.1) yields
$ex(p;K_{1,3})=p$. So the result is true for $n=4$.
\par Now we assume $n\geq 5$. By (2.1), it suffices to show that
$ex(p;K_{1,n-1})\ge [\f{(n-2)p}2]$. Set $k=[\f{p+1}2]$,
$V=\{1,2,\ldots,2k\}$ and $M=\{12,34,\cdots,(2k-1)(2k)\}$.
 Let us
consider the following four cases.
\par {\bf Case 1.} $2\mid p$ and $2\nmid n$. Set $G=(V,E)$, where
$$E=\big\{ij\ |\ i,j\in V,\ j-i\in\{1,2k-1,k,k\pm 1,\ldots,k\pm
(n-5)/2\}\big\}.$$
 Clearly
$G$ is an $(n-2)$-regular graph of order $p$ and so $G$ does not
contain $K_{1,n-1}$. Hence, $ex(p;K_{1,n-1})\ge e(G)
=\f{(n-2)p}2=[\f{(n-2)p}2]$.
 \par {\bf Case 2.} $2\mid p$ and $2\mid n$. Set
$$E_1=\big\{ij\ |\ i,j\in V,\
j-i\in\{1,2k-1,k,k\pm 1,\ldots,k\pm (n-4)/2\}\big\}.$$ Then
$M\subset E_1$. Let $G=(V,E_1-M)$. We see that $G$ is an
$(n-2)$-regular graph of order $p$ and so $G$ does not contain
$K_{1,n-1}$. Hence, $ex(p;K_{1,n-1})\ge e(G)
=\f{(n-2)p}2=[\f{(n-2)p}2]$.
\par {\bf Case 3.} $2\nmid p$ and $2\mid n$.
Let $G$ be the $(n-2)$-regular graph of order $2k$ constructed in
Case 2. Let
$$v_1=k-\f n2+3,\ v_2=k-\f n2+4,\ldots,v_{n-3}=k+\f n2-1\qtq{and}v_{n-2}=2k.$$
Then clearly $v_1,\ldots,v_{n-2}$ are all the vertices adjacent to
the vertex $1$. If $2\mid k-\f n2$, then $v_1,v_3,\ldots,v_{n-5}$
are odd and so $v_1v_2,v_3v_4,\ldots,v_{n-5}v_{n-4}\in M$. Thus,
$v_1v_2,v_3v_4,\ldots,$ $v_{n-5}v_{n-4}\notin E(G)$. As $2k-(k+\f
n2-1)=k-\f{n-2}2$, we see that $v_{n-3}v_{n-2}\not\in E_1$ and so
$v_{n-3}v_{n-2}\not\in E(G)$.
 Let
 $$G'=G-\{1\}+\{v_1v_2,v_3v_4,\ldots,v_{n-5}v_{n-4},v_{n-3}v_{n-2}\}.$$
We see that $G'$ is an $(n-2)$-regular graph of order $p$. Hence,
$ex(p;K_{1,n-1})\ge e(G')=\f{(n-2)p}2=[\f{(n-2)p}2]$. \par If
$2\nmid k-\f n2$, then $v_2,v_4,\ldots,v_{n-4}$ are odd and so
$v_2v_3,v_4v_5,\ldots,v_{n-4}v_{n-3}\in M$. Thus,
$v_2v_3,v_4v_5,\ldots,v_{n-4}v_{n-3}\notin E(G)$. As $p+1=2k>n$ we
have $k-\f n2+3>3$ and so $2,3\notin\{v_1,\ldots,v_{n-2}\}$. Clearly
$2v_{n-2},3v_1\notin E_1$ and so $2v_{n-2},3v_1\notin E(G)$. Let
$$G'=G-\{1\}-\{23\}+\{v_2v_3,v_4v_5,\ldots,v_{n-4}v_{n-3},3v_1,2v_{n-2}\}.$$
Then $G'$ is an $(n-2)$-regular graph of order $p$. Hence,
$ex(p;K_{1,n-1})\ge e(G')=\f{(n-2)p}2=[\f{(n-2)p}2]$.
\par {\bf Case 4.} $2\nmid p$ and $ 2\nmid n$. As $2\mid n+1$,
we can construct an  $(n-1)$-regular graph $G_1$ of order $p$ by
using the argument in Case 3. Let
$$M_1=\cases \{23,45,\ldots,(2k-2)(2k-1),k(2k)\}&\t{if $2\mid
k-\f{n+1}2$,}\\
\{2(2k),3(k+3-\f{n+1}2),45,67,\ldots,(2k-2)(2k-1)\} &\t{if $2\nmid
k-\f{n+1}2$.}\endcases $$
 It is easily seen that $M_1\subset G_1$. Set $G_2=G_1-M_1$.
 Then for $i=2,3,\ldots,2k$ we
have $$d_{G_2}(i)=\cases n-3&\t{if $2\mid k-\f{n+1}2$ and $i=k$,
or if $2\nmid k-\f{n+1}2$ and $i=k+3-\f{n+1}2$,}
\\n-2&\t{otherwise.}\endcases$$
Thus $G_2$ does not contain $K_{1,n-1}$ and
$$2e(G_2)=\sum_{i=2}^{2k}d_{G_2}(i)=n-3+(2k-2)(n-2)=(n-2)p-1.$$
Hence $ex(p;K_{1,n-1})\ge e(G_2)=\f{(n-2)p-1}2=[\f{(n-2)p}2].$
\par Putting all the above together we prove the theorem.
\end{proof}
\pro{\noindent Corollary 2.1} Let $k,p\in\Bbb N$ with $p\ge k+2$. Then there
exists a $k-$regular graph of order $p$ if and only if $2\mid
kp$.
\endpro
%   Proof.
\begin{proof}
If $G$ is a $k-$regular graph of order $p$, then
$kp=2e(G)$ and so $2\mid kp$. If $2\mid kp$, by the proof of Theorem
2.1 we know that there exists a $k-$regular graph of order $p$.
\end{proof}
\par{\it \noindent Remark $2.1.$} In [4] Kirkman showed that $K_{2n}$ is 1-factorable.
In [5] Petersen proved that a graph $G$ is 2-factorable if and only
if $G$ is $2p$-regular. Thus, Corollary 2.1 can be deduced from [4]
and [5].

 \section{The Evaluation of $ex(p;T_n)$}
\pro{\noindent Theorem 3.1} Let $p,n\in\Bbb N$ with $p\geq n\geq 5$. Let
$r\in\{ 0,1,\ldots,n-2\}$ be given by $p\e r\mod{n-1}$. Then
$$ex(p;T_n)= \cases \big[\f{(n-2)(p-1)-r-1}2\big] &\t{if}\
n\ge 7\ \t{and}\
2\le r\le n-4,\\
\f{(n-2)p-r(n-1-r)}2&\t{otherwise.}\endcases$$
\endpro
%   Proof.
\noindent{\it Proof.\/} Let $G$ be an extremal graph of order $p$
not containing $T_n$. Suppose $v_0\in V(G)$ and $G_0$ is the
component of $G$ such that $v_0\in V(G_0)$. If $d(v_0)=m\geq n-1$,
as $G$ does not contain $T_n$ we see that $G_0$ is a copy of
$K_{1,m}$. Suppose $m+1=k'(n-1)+r'$ with $k'\in\Bbb N$ and
$r'\in\{0,1,\ldots,n-2\}.$
 Then $k'K_{n-1}\cup K_{r'}$ does not contain $T_n$.
  As $\f{n-2}2>1$ and $\binom {r'}2-(r'-1)=\f{(r'-1)(r'-2)}2\ge 0$, we find
$$ e(k'K_{n-1}\cup K_{r'})=k'\binom{n-1}2+\binom{r'}2>k'(n-1)+r'-1=m=e(K_{1,m})=e(G_0).$$
Hence $G_0\notin Ex(m+1;T_n)$ and so $G\notin Ex(p;T_n)$. This
contradicts the assumption. Therefore
 $d(v_0)\leq n-2$ and so $\Delta(G)\le n-2$. If $d(v_0)=n-2$, as
$G_0$ is an extremal graph not containing $T_n$ we see that $G_0$ is
a copy of $K_{n-1}$.
\par Suppose $p=k(n-1)+r$. Then $k\in\Bbb N$. From the above we may
assume $G=sK_{n-1}\cup G_1$ with $s\in\{0,1,\ldots,k\}$ and
$\Delta(G_1)\le n-3$. If $s=k$, then clearly $G_1=K_r$ and so
$e(G)=k\b{n-1}2+\b r2$. If $s\le k-1$, as $\Delta(G_1)\le n-3$
implies $G_1$ does not contain any copies of $T_n$, we see that
$G_1\in Ex((k-s)(n-1)+r;K_{1,n-2})$. By Theorem 2.1 we have
$e(G_1)=[\f{(n-3)((k-s)(n-1)+r)}2]$. Hence
$$e(G)=e(sK_{n-1}\cup
G_1)=s\b{n-1}2+\Big[\f{(n-3)((k-s)(n-1)+r)}2\Big].$$
   Set $f(x)=x\binom{n-1}2+[\f{(n-3)((k-x)(n-1)+r)}2]$.
  Then
   $$\aligned f(x+1)&=(x+1)\binom{n-1}2+\[\f{(n-3)((k-x)(n-1)+r)-(n-3)(n-1)}2\]
  \\&=x\binom{n-1}2+\[\f{(n-3)((k-x)(n-1)+r)}2+\f{n-1}2\]
  >f(x).\endaligned$$ Thus, $f(k-1)>f(k-2)>\ldots>f(0)$.
  Since $G$ is an extremal graph, by the above we must have
  $s=k-1$ or $k$ and so
  $$\align&ex(p; T_n)=e(G)\\&=\max\bigg\{(k-1)\binom{n-1}2+\[\f{(n-3)(n-1+r)}2\],k\binom{n-1}2+\binom
  r2\bigg\}.\endalign$$
Observe that
$$\f{(n-3)(n-1+r)}2-\f{r(r-1)}2-\f{(n-1)(n-2)}2=\f{r(n-2-r)-(n-1)}2.$$
 We then have
$$ex(p;T_n)=k\binom{n-1}2+\binom
r2+\max\bigg\{0,\Big[\f{r(n-2-r)-(n-1)}2\Big]\bigg\}.$$ If
$r\in\{1,n-3,n-2\}$, then clearly $[\f{r(n-2-r)-(n-1)}2]<0.$ For
$n=6$ and $r=2$ we also have $[\f{r(n-2-r)-(n-1)}2]=-1<0.$ Now
assume $n\geq 7$ and $2\leq r\leq n-4$ . Then
$$\aligned r(n-2-r)-(n-1)&=\f{n^2-8n+8}4-\Big(r-\f{n-2}2\Big)^2
\\&\geq \f{n^2-8n+8}4-\Big(2-\f{n-2}2\Big)^2=n-7\ge 0\endaligned$$
 and so $[\f{r(n-2-r)-(n-1)}2]\ge 0.$ Hence
$$\aligned &ex(p;T_n)\\&=\cases k\b{n-1}2+\b r2+\big[\f{r(n-2-r)-(n-1)}2\big]&\t{if $n\ge
7$ and $2\le r\le n-4$,}\\k\b{n-1}2+\b r2&\t{otherwise.}
\endcases\endaligned$$
To see the result, we note that $k\b{n-1}2+\b
r2=\f{(n-2)(p-r)+r^2-r}2=\f{(n-2)p-r(n-1-r)}2$ and
$$k\b{n-1}2+\b r2+\Big[\f{r(n-2-r)-(n-1)}2\Big]
=\Big[\f{(n-2)(p-1)-r-1}2\Big].\hskip+2cm \square $$

 \section{The Evaluation of $ex(p;T_n^*)$}
 \par For $n\ge 4$ we recall that $T_n^*=(V,E)$ is the tree on $n$ vertices with
 $V=\{v_0,v_1,\ldots,v_{n-1}\}$ and
 $E=\{v_0v_1,\ldots,v_0v_{n-3},v_{n-3}v_{n-2},v_{n-2}v_{n-1}\}$.
Clearly $T_4^*=P_4$ and $T_5^*=P_5$.
 \pro{\noindent Lemma 4.1} Let $p,n\in\Bbb
N$ with $p\ge n\ge 6$, and let $G\in Ex(p;T_n^*)$. Then $\Delta
(G)\le n-2.$
\endpro\begin{proof}
Suppose that $v_0\in V(G), d(v_0)=m\ge n-1$ and $\Gamma
(v_0)=\{v_1,\ldots,v_m\}.$  Let $G_0$ be the component of $G$ with
$v_0\in V(G_0).$ If there are exactly $t$ vertices
$u_1,\ldots,u_t\in V(G_0)$ such that
$d(u_1,v_0)=\cdots=d(u_t,v_0)=2,$ then clearly
$d(u_1)=\cdots=d(u_t)=1$,
$V(G_0)=\{v_0,v_1,\ldots,v_m,u_1,\ldots,u_t\}$ and $|V(G_0)|=1+m+t$.
If $u_iv_j\notin E(G_0)$ for some $j\in\{1,2,\ldots,m\}$ and every
$i=1,2,\ldots,t$, then clearly $d(v_j)\le 2$. Thus, $e(G_0)\le
m+t+\f m2$. Set $1+m+t=k(n-1)+r(0\le r<n-1).$ We see that
$$\aligned &k\b{n-1}2+\b r2-\f{3m}2-t\\&=
\f{(n-2)(1+m+t-r)+r(r-1)-3m-2t}2\\&=\f{(m+t)(n-5)-r(n-1-r)+(n-2)+t}2\\&\ge
\f{(n-1)(n-5)+(n-2)-r(n-1-r)}2\\&\ge\f{(n-1)(n-5)+n-2-\f{(n-1)^2}4}2
=\f{3(n-3)^2-16}8
>0.\endaligned$$ Since $kK_{n-1}\cup K_r$ does not contain any
copies of $T_n^*,$ applying the above we deduce
$$e(G_0)\le\f{3m+2t}2<k\b{n-1}2+\b r2=e(kK_{n-1}\cup K_r)\le
ex(1+m+t;T_n^*).$$ As $G$ is an extremal graph not containing
$T_n^*,$ we must have $e(G_0)=ex(1+m+t; T_n^*)$. This contradicts
the above inequality $e(G_0)<ex(1+m+t;T_n^*)$. Hence the assumption
$d(v_0)\ge n-1$ is not true. Thus $\Delta(G)\le n-2.$ The proof is
now complete.
\end{proof}
 \pro{\noindent Lemma 4.2} Let
$p,n\in\Bbb N$ with $p\ge n\ge 5$, and let $G\in Ex(p;T_n^*)$.
Suppose that $v_0\in V(G), d(v_0)=n-2$ and $G_0$ is the component of
$G$ such that $v_0\in V(G_0)$. Then $G_0\cong K_{n-1}.$
\endpro
%   Proof.
\begin{proof}
Suppose $\Gamma(v_0)=\{v_1,\ldots,v_{n-2}\}$ and there are
exactly $t$ vertices $u_1,\ldots,u_t\in V(G_0)$ such that
$d(u_1,v_0)=\cdots=d(u_t,v_0)=2.$ We first assume $t\ge 1$. Then
clearly $d(u_1)=\cdots=d(u_t)=1$ and
$V(G_0)=\{v_0,v_1,\ldots,v_{n-2},u_1,\ldots,u_t\}.$ If $u_1v_i\in
E(G)$ for some $i\in\{1,2,\ldots,n-2\}$, then clearly $v_iv_j\notin
E(G)$ for all $j\in\{1,2,\ldots,n-2\}\setminus\{i\}.$ Thus,
$$e(G_0)\le n-2+t+\b{n-2-t}2\le \b{n-2}2+t+1.$$ Assume $t=q(n-1)+t_0$
with $q\in\Bbb Z$ and $t_0\in\{0,1,\ldots,n-2\}.$ Then $$\aligned
&e((1+q)K_{n-1}\cup
K_{t_0})-\b{n-2}2-t-1\\&=(1+q)\b{n-1}2+\b{t_0}2-\b{n-2}2-q(n-1)-t_0-1
\\&=n-4+q\f{(n-1)(n-4)}2+\f{(t_0-1)(t_0-2)}2>0.\endaligned$$ As $(1+q)K_{n-1}\cup
K_{t_0}$ does not contain $T_n^*,$ applying the above we get
$$\aligned e(G_0)\le\b{n-2}2+t+1<e((1+q)K_{n-1}\cup K_{t_0})\le
ex(n-1+t;T_n^*).\endaligned$$ Since $G_0$ is an extremal graph of
order $n-1+t$ not containing $T_n^*,$ we must have
$e(G_0)=ex(n-1+t;T_n^*).$ This contradicts the above assertion. So
$t\ge 1$ is not true and hence $V(G_0)=\{v_0,v_1,\ldots,v_{n-2}\}.$
As $G_0$ is an extremal graph not containing $T_n^*,$ we see that
$G_0\cong K_{n-1}.$ This proves the lemma.
\end{proof}
\pro{\noindent Lemma 4.3} Let
$n,t\in\Bbb N$ with $n\ge 4$, and let $G\in Ex(n-2+t;T_n^*)$.
Suppose that $G$ is connected and $\Delta(G)=n-3.$
 Then $t\le n-4$ and $e(G)\le (n-3)^2.$
 \endpro
 %  Proof.
 \begin{proof}
 Suppose $v_0\in V(G),d(v_0)=n-3,\Gamma(v_0)=\{v_1,\ldots,v_{n-3}\}$ and
 $V(G)=\{v_0,v_1,\ldots,$ $v_{n-3},u_1,\ldots,u_t\}$. Then $d(u_i,v_0)=2$ and $u_1,\ldots,u_t$
  must be independent. As $G$ is connected and $u_i$ is adjacent to
 some vertex in $\Gamma(v_0),$
  we have
  $$e(G)\le\sum_{i=1}^{n-3}d(v_i)\le\sum_{i=1}^{n-3}(n-3)=(n-3)^2.$$
  On the other hand,
  $$e(K_{n-1}\cup K_{n-4})=\f{(n-1)(n-2)+(n-4)(n-5)}2=n^2-6n+11>(n-3)^2.$$
  Thus, for $t\ge n-3$ we have
  $$\aligned e(G)&=ex(n-2+t;T_n^*)\ge e\(K_{n-1}\cup K_{n-4}\cup (t-(n-3))K_1\)
  \\&=e(K_{n-1}\cup K_{n-4})>(n-3)^2.\endaligned$$
  This contradicts the fact $e(G)\le (n-3)^2.$ So $t\le n-4.$ The proof is now complete.
\end{proof}
\pro{\noindent Lemma 4.4} Let $p,n\in\Bbb N$ with $p\ge n\ge 4,$ and $G\in
Ex(p;T_n^*)$. Suppose $\Delta(G)\le n-3$. Then $p\le 2n-6.$
\endpro
%   Proof.
\begin{proof}
Assume $p=2n-4+t.$ If $t\ge 2n,$ we may write $t-2=k(n-1)+r,$
where $k\in\Bbb N$ and $r\in\{0,1,\ldots,n-2\}.$ Let $G_0\in
Ex(n-1+r;K_{1,n-3}).$ From Theorem 2.1 we have
$e(G_0)=[\f{(n-1+r)(n-4)}2]$. Clearly $k(n-1)=t-2-r\ge 2n-2-r>r+1.$
Thus,
$$\aligned e((k+1)K_{n-1}\cup
G_0)&=(k+1)\b{n-1}2+\[\f{(n-1+r)(n-4)}2\]
\\&\ge\f{(k+1)(n-1)(n-2)+(n-1+r)(n-4)-1}2
\\&=\f{\((k+2)(n-1)+r\)(n-3)}2+\f{k(n-1)-r-1}2
\\&>\f{\((k+2)(n-1)+r\)(n-3)}2=\f{(n-3)p}2.\endaligned$$
On the other hand, as $(k+1)K_{n-1}\cup G_0$ does not contain
$T_n^*,$ we have $$e((k+1)K_{n-1}\cup G_0)\le
ex(p;T_n^*)=e(G)\le\f{(n-3)p}2.$$ This is a contradiction. Hence
$t<2n.$
\par If $t=2n-1,$ then $p=2n-4+t=3(n-1)+n-2$ and so
$$\f{(n-3)p}2<e(3K_{n-1}\cup K_{n-2})\le
ex(p;T_n^*)=e(G)\le\f{(n-3)p}2.$$ This is also a contradiction.
\par If $n-1\le t<2n-1,$ setting $G_0\in Ex(t-2;K_{1,n-3})$ and
using Theorem 2.1 we see that
$$e(G_0)=ex(t-2;K_{1,n-3})=\[\f{(n-4)(t-2)}2\].$$
It is clear that $2K_{n-1}\cup G_0$ does not contain $T_n^*$ as a
subgraph and
$$\aligned e(2K_{n-1}\cup G_0)&=2\b{n-1}2+\[\f{(n-4)(t-2)}2\]
\\&\ge (n-1)(n-2)+\f{(n-4)(t-2)-1}2
\\&=\f{(2n-4+t)(n-3)}2+\f{2n-1-t}2
>\f{(2n-4+t)(n-3)}2.\endaligned$$ On the other hand,
$$e(2K_{n-1}\cup G_0)\le
ex(2n-4+t;T_n^*)=e(G)\le\f{(2n-4+t)(n-3)}2.$$ This is a
contradiction.
\par By the above, we may assume $t\le n-2.$ If $t=n-2$, then
$$\aligned ex(3n-6;T_n^*)&\ge e(2K_{n-1}\cup
K_{n-4})=\f{2(n-1)(n-2)+(n-4)(n-5)}2\\&>\f{(3n-6)(n-3)}2\ge
e(G)=ex(3n-6;T_n^*).\endaligned$$ This is a contradiction. If
$t=n-3$, then $$ \aligned ex(3n-7;T_n^*)&\ge e(K_{n-1}\cup
K_{n-3,n-3})=\f{(n-1)(n-2)}2+(n-3)^2\\&>\f{(3n-7)(n-3)}2\ge
e(G)=ex(3n-7;T_n^*).\endaligned$$ This is also a contradiction. Thus
$t\not=n-2,n-3$.
\par Now we
assume that $1\le t\le n-4.$ Suppose $H\in Ex(n-3;K_{1,n-3-t})$ and
$V(H)=\{v_1,\ldots,v_{n-3}\}.$ We construct a graph
$G_0=(V(G_0),E(G_0))$ of order $n-3+t$ by defining
$V(G_0)=\{u_1,\ldots,u_t\}\cup V(H)$ and $E(G_0)=\{u_iv_j:1\le i\le
t,1\le j\le n-3\}\cup E(H).$ It is easily seen that $d_{G_0}(v_i)\le
n-4(1\le i\le n-3)$ and so $G_0$ does not contain any copies of
$T_n^*.$ Hence, $$\align e(K_{n-1}\cup G_0)&=\b{n-1}2+e(G_0)\\&\le
ex(2n-4+t;T_n^*)=e(G)\le\f{(2n-4+t)(n-3)}2.\endalign$$ Using Theorem
2.1 we see that
$$\aligned
e(G_0)&=(n-3)t+\[\f{(n-3)(n-4-t)}2\]\\&\ge(n-3)t+\f{(n-3)(n-4-t)-1}2
\\&=\f{(2n-4+t)(n-3)}2-\b{n-1}2+\f
12\\&>\f{(2n-4+t)(n-3)}2-\b{n-1}2,\endaligned$$ this contradicts the
above assertion.
\par By the above we have $t\le 0$ and so $p\le 2n-4.$ If $p=2n-4$,
since $K_{n-1}\cup K_{n-3}$ does not contain $T_n^*$ we have
$$\aligned ex(2n-4;T_n^*)&\ge e(K_{n-1}\cup
K_{n-3})=\f{(n-1)(n-2)+(n-3)(n-4)}2\\&>\f{(2n-4)(n-3)}2\ge
e(G)=ex(2n-4;T_n^*).\endaligned$$ This is a contradiction.
\par Now we assume $p=2n-5.$ It is clear that
$$e(K_{n-1}\cup K_{n-4})=\f{(n-1)(n-2)+(n-4)(n-5)}2=n^2-6n+11.$$
 As $K_{n-1}\cup K_{n-4}$ does not contain $T_n^*,$ we see that $n^2-6n+11
 \le ex(2n-5;T_n^*)=e(G).$ If $\Delta(G)\le n-4,$ then clearly
 $e(G)\le\f{(2n-5)(n-4)}2<n^2-6n+11$. This is a contradiction.
 Hence, $\Delta(G)=n-3.$ Suppose that $G_1$ is the component of
 $G$ such that $\Delta(G_1)=n-3.$ If $|V(G_1)|=n-2+s$ for some
 $s\in\{0,1,\ldots,n-3\},$ by Lemma 4.3 we have $s\le n-4.$ As $G$
 is an extremal graph we have $G\backslash G_1\cong K_{n-3-s}$ and
 so $$\aligned e(G)&=e(G_1)+e(G\backslash G_1)\le\f{(n-2+s)(n-3)}2+\b{n-3-s}2
 \\&=\f 12\Big(s-\f{n-4}2\Big)^2+\f{7n^2-40n+56}8
 \\&\le\f
 12\Big(\f{n-4}2\Big)^2+\f{7n^2-40n+56}8=n^2-6n+9<n^2-6n+11,\endaligned$$
 this contradicts the above assertion $e(G)\ge n^2-6n+11.$ Therefore
 $p\neq 2n-5$ and so $p\le 2n-6$, which completes the proof.
\end{proof}
\pro{\noindent Theorem 4.1} Let $p,n\in\Bbb N$ with $p\ge n-1\ge 5$, and let
$p=k(n-1)+r$ with $k\in\Bbb N$ and $r\in\{0,1,\ldots,n-2\}.$ Then
$$\aligned
&ex(p;T_n^*)\\&=\cases\f{(k-1)(n-1)(n-2)}2+ex(n-1+r;T_n^*)&\t{if}\
1\le r\le n-5;\\\f{(n-2)p-r(n-1-r)}2&\t{if}\
r\in\{0,n-4,n-3,n-2\}.\endcases\endaligned$$\endpro
\begin{proof}
Suppose
$m\in\Bbb N$ and $m\ge 2n-5$. We assert that
$$ex(m;T_n^*)=\f{(n-1)(n-2)}2+ex(m-(n-1);T_n^*).\tag 4.1$$ Assume $G\in
Ex(m;T_n^*).$ From Lemma 4.1 we know that $\Delta(G)\le n-2.$ As
$m\ge 2n-5,$ by Lemma 4.4 we have $\Delta(G)=n-2.$ Using Lemma 4.2
we see that $G$ has a component isomorphic to $K_{n-1}$ and so (4.1)
is true. From (4.1) we deduce that for $k\ge 2$,
$$ \aligned &ex(p;T_n^*)-ex(n-1+r;T_n^*)
\\&=\sum_{s=1}^{k-1}\big\{ex((s+1)(n-1)+r;T_n^*)-ex(s(n-1)+r;T_n^*)\big\}
=(k-1)\b{n-1}2.\endaligned$$ This is also true for $k=1$. \par For
$r=0,$ we have $ex(n-1+r;T_n^*)=e(K_{n-1})=\b{n-1}2$ and so
$$ex(p;T_n^*)=(k-1)\b{n-1}2+\b{n-1}2=k\b{n-1}2=\f{(n-2)p}2.$$ For
$r\in\{n-4,n-3,n-2\}$ we have $n-1+r\ge 2n-5$ and so by (4.1)
$$\aligned ex(p;T_n^*)&=(k-1)\b{n-1}2+ex(n-1+r;T_n^*)
\\&=(k-1)\b{n-1}2+\b{n-1}2+ex(r;T_n^*)
=k\b{n-1}2+e(K_r)\\&=\f{(n-2)(p-r)}2+\b r2
=\f{(n-2)p-r(n-1-r)}2\endaligned$$ as asserted. The proof is now
complete.
\end{proof}
 \pro{\noindent Theorem 4.2} Let $p,n\in\Bbb N$ with
$p\ge n\ge 6$ and $p=k(n-1)+1$ with $k\in\Bbb N$. Then
$$ex(p;T_n^*)=\f{(n-2)(p-1)}2.$$\endpro
%   Proof.
\begin{proof}
 Let $G_0\in Ex(n;T_n^*).$ If $\Delta(G_0)\le n-3,$ then
$e(G_0)\le\f{(n-3)n}2<\f{(n-1)(n-2)}2$. On the other hand,
$e(G_0)=ex(n;T_n^*)\ge e(K_{n-1}\cup K_1)=\f{(n-1)(n-2)}2$. This is
a contradiction. Thus $\Delta(G_0)\ge n-2$. Applying Lemmas 4.1 and
4.2 we see that $G_0\cong K_{n-1}\cup K_1$ and so
$ex(n;T_n^*)=e(G_0)=\f{(n-1)(n-2)}2.$ Now applying Theorem 4.1 we
obtain
$$ex(p;T_n^*)=\f{(k-1)(n-1)(n-2)}2+ex(n;T_n^*)=k\b{n-1}2=\f{(n-2)(p-1)}2.$$
This is the result.
\end{proof}
 \pro{\noindent Theorem 4.3} Let $p,n\in\Bbb N$, $p\ge n\ge
7$ and $p=k(n-1)+n-5$ with $k\in\Bbb N$. Then
$$ex(p;T_n^*)=\f{(n-2)(p-2)}2+1.$$
\endpro
%   Proof.
\begin{proof}
Let $G_0\in Ex(2n-6;T_n^*).$ If $\Delta(G_0)\le n-3,$ then
$e(G_0)\le\f{(n-3)(2n-6)}2=(n-3)^2$. As $K_{n-3,n-3}$ does not
contain any copies of $T_n^*$, we see that $e(G_0)\ge
e(K_{n-3,n-3})=(n-3)^2$. Hence $e(G_0)=(n-3)^2$. If $\Delta(G_0)\ge
n-2$, by Lemmas 4.1 and 4.2 we have $G_0\cong K_{n-1}\cup K_{n-5}$
Thus, $e(G_0)=e(K_{n-1}\cup K_{n-5})=\b {n-1}2+\b{n-5}2=n^2-7n+16$.
Since $(n-3)^2=n^2-6n+9\ge n^2-7n+16$, we see that
$ex(2n-6;T_n^*)=(n-3)^2$. Now applying the above and Theorem 4.1 we
deduce
$$\align ex(p;T_n^*)&=(k-1)\b{n-1}2+ex(2n-6;T_n^*)=(k-1)\b{n-1}2+(n-3)^2
\\&=k\f{(n-1)(n-2)}2+\f{n^2-9n+16}2=\f{(n-2)(p-2)}2+1.\endalign$$
 This is the result.
\end{proof}
 \pro{\noindent Lemma 4.5} Let $n,r\in\Bbb N$ with $n\ge 7$ and
$r\le n-5$. Then there is an extremal graph $G\in
Ex(n-1+r;\{K_{1,n-2},T_n^*\})$ such that $\Delta(G)=n-3$ and $G$ is
connected.
 \endpro
%    Proof.
\begin{proof}
Let $G\in
Ex(n-1+r;\{K_{1,n-2},T_n^*\})$. Then $\Delta(G)\le n-3$. For $r=n-5$
we see that $K_{n-3,n-3}\in Ex(n-1+r;\{K_{1,n-2},T_n^*\})$. So the
result is true.
\par Now we assume $r\le n-6.$
Suppose $H\in Ex(n-3;K_{1,n-5-r})$ and
$V(H)=\{v_1,\ldots,v_{n-3}\}.$ From Theorem 2.1 we know that
$e(H)=ex(n-3;K_{1,n-5-r})=[\f{(n-3)(n-6-r)}2]$. Now we construct a
graph $G_0=(V(G_0),E(G_0))$ of order $n-1+r$ by defining
$V(G_0)=\{u_0,\ldots,u_{r+1}\}\cup V(H)$ and $E(G_0)=\{u_iv_j:\ 0\le
i\le r+1,1\le j\le n-3\}\cup E(H).$ It is easily seen that
 $d_{G_0}(v_i)\le n-4(1\le i\le n-3)$, $\Delta(G_0)=n-3$ and so $G_0$
does not contain any copies of $T_n^*$ and $K_{1,n-2}$.
 Thus, for any $G\in
Ex(n-1+r;\{K_{1,n-2},T_n^*\})$,
$$e(G)\ge
e(G_0)=(n-3)(r+2)+\[\f{(n-3)(n-6-r)}2\]=\[\f{(n-3)(n-2+r)}2\].$$ If
$\Delta(G)\le n-4$, we must have  $G\in Ex(n-1+r;K_{1,n-3})$ and so
$e(G)=[\f{(n-4)(n-1+r)}2]$ by Theorem 2.1. As $G$ is an extremal
graph and
$$\align \[\f{(n-3)(n-2+r)}2\]&\ge \f{(n-3)(n-2+r)-1}2=\f{(n-4)(n-1+r)+r+1}2
\\& >\f{(n-4)(n-1+r)}2\ge \[\f{(n-4)(n-1+r)}2\],\endalign$$
by the above we must have $\Delta(G)=n-3$.
\par Now assume $\Delta(G)=n-3$. If $G$ is connected, the result is
true. Suppose that $G$ is not connected. Let $G_1$ be a component of
$G$ with $\Delta(G_1)=n-3$ and $|V(G_1)|=n-1+r-s.$ Then $1\le s\le
r+1\le n-5$. As $G$ is an extremal graph, we must have $G=G_1\cup
K_s$. Thus,
$$e(G)=e(G_1)+\b s2\le\Big[\f{(n-3)(n-1+r-s)}2\Big]+\f{s(s-1)}2.$$
On the other hand, $e(G)\ge e(G_0)=[\f{(n-3)(n-2+r)}2].$ Therefore,
$$\Big[\f{(n-3)(n-2+r)}2\Big]-\Big[\f{(n-3)(n-1+r-s)}2\Big]-\f{s(s-1)}2\le
0.$$
 For $s\ge 2$ we have
$(s-1)(n-3-s)=(s-2)(n-4-s)+n-5\ge n-5$ and so
$$\align &\Big[\f{(n-3)(n-2+r)}2\Big]-\Big[\f{(n-3)(n-1+r-s)}2\Big]-\f{s(s-1)}2
\\&\ge \Big[-\f{s^2-(n-2)s+n-3}2\Big]=\Big[\f{(s-1)(n-3-s)}2\Big]\ge \Big[\f{n-5}2\Big]>0.
\endalign
$$ This contradicts the previous inequality. Thus $s=1$ and hence
$e(G)=e(G_1)\le [\f{(n-3)(n-2+r)}2]=e(G_0)$. By the previous
argument, $e(G)\ge e(G_0)$.
 Therefore $e(G)=e(G_0)$. As $G_0$ is
connected and $\Delta(G_0)=n-3$, we see that the result is true.
\end{proof}
\pro{\noindent Lemma 4.6} Let $n,r\in\Bbb N$ with $n\ge 11$ and $3\le r\le
n-5$.  Then there is an extremal graph $G\in Ex(n-1+r;T_n^*)$ such
that $\Delta(G)=n-3$ and $G$ is connected. Moreover,
$ex(n-1+r;T_n^*)=ex(n-1+r;\{K_{1,n-2},T_n^*\})$.
\endpro
%   Proof.
\begin{proof}
Let $G\in Ex(n-1+r;T_n^*)$. For $r=n-5$ let
$G_0=K_{n-3,n-3}$. For $r\le n-6$ let $G_0$ be the graph constructed
in the proof of Lemma 4.5. Then $\Delta(G_0)=n-3$ and $G_0$ does not
contain any copies of $T_n^*$. Thus, $e(G)\ge e(G_0)$. For $r=n-5$
we have $e(G_0)=(n-3)^2$. For $r\le n-6$ we have
$e(G_0)=[\f{(n-3)(n-2+r)}2]$. Since $(n-3)^2\ge
\f{(n-3)(n-2+n-5)}2$, we always have $e(G)\ge [\f{(n-3)(n-2+r)}2]$
for $r\le n-5$.
\par If $\Delta(G)\ge n-2$, by Lemmas 4.1 and 4.2 we have $G\cong
K_{ n-1}\cup K_r$. Thus, $e(G)=\b{n-1}2+\b r2$. Since $3\le r\le
n-5$ and $n\ge 11$ we see that  $(r-2)(n-4-r)\ge 4$ and so
$$\[\f{(n-3)(n-2+r)}2\]-\b{n-1}2-\b
r2=\[\f{(r-2)(n-4-r)-4}2\]\ge 0.$$ Therefore $e(G)\le e(G_0)$ and so
$e(G)=e(G_0)$. Since $\Delta(G_0)=n-3$ and $G_0$ is connected, the
result holds in this case.
\par Now we assume $\Delta(G)\le n-3$. Then $G\in Ex(n-1+r;\{K_{1,n-2},T_n^*\})$. Applying
Lemma 4.5 we see that the result is true. Thus the lemma is proved.
\end{proof}
 \pro{\noindent Lemma 4.7} Let
$n,r\in\Bbb N$ with $n\ge 7$ and $r\le n-5$.  Then
$$ex(n-1+r;\{K_{1,n-2},T_n^*\})=(n-3)(r+2)+ex(n-3;\{K_{1,n-4-r},T_{n-2-r}^*\}).$$
Moreover, for $r\ge \f{n-7}2$ we have
$$\align&ex(n-1+r;\{K_{1,n-2},T_n^*\})
\\&=(n-3)(r+2)+\max\big\{(n-5-r)^2,\big[\f{(n-6-r)(n-3)
}2\big]\big\}.\endalign$$\endpro
%    Proof.
\begin{proof}
It is clear that
 $ex(2n-6;\{K_{1,n-2},T_n^*\})=e(K_{n-3,n-3})=(n-3)^2$. So the
 result is true for $r=n-5$.
\par Now assume $r\le n-6$. By Lemma 4.5, we can choose a graph
$G\in Ex(n-1+r;\{K_{1,n-2},T_n^*\})$ so that $ \Delta(G)=n-3$ and
$G$ is connected. Suppose $u_0\in
V(G),d(u_0)=n-3,\Gamma(u_0)=\{v_1,\ldots,v_{n-3}\}$ and
$V(G)=\{v_1,\ldots,v_{n-3},u_0,u_1,\ldots,u_{r+1}\}.$ Then
$d(u_i,u_0)=2$ for $i=1,2,\ldots,r+1$ and
$\{u_0,u_1,\ldots,u_{r+1}\}$ is an independent set.
 If $u_iv_j\notin
E(G)$ for some $i\in\{1,2,\ldots,r+1\}$ and
$j\in\{1,2,\ldots,n-3\}$, as $G$ is an extremal graph we see that
$v_jv_k\in E(G)$ for some $k\in\{1,2,\ldots,n-3\}-\{j\}$. Set
$G_1=G-v_jv_k+u_iv_j$. Then clearly $G_1$ does not contain $T_n^*$,
$e(G)=e(G_1)$, $\Delta(G_1)=n-3$ and $G_1$ is connected. Repeating
the above step we see that there is an extremal graph $G'\in
Ex(n-1+r;\{K_{1,n-2},T_n^*\})$ such that
$V(G')=\{v_1,\ldots,v_{n-3},u_0,u_1,\ldots,u_{r+1}\}$,
$\Gamma(u_i)=\{v_1,\ldots,v_{n-3}\}$ for $i=0,1,\ldots,r+1$,
$\Delta(G')=n-3$ and $G'$ is connected. It is easily seen that
$$e(G')=(n-3)(r+2)+e(G'[v_1,\ldots,v_{n-3}]).$$
 Set $H=G'[v_1,\ldots,v_{n-3}].$
Since $\Delta(G')=n-3$ and $G'\in Ex(n-1+r;\{K_{1,n-2},T_n^*\})$, we
see that $\Delta(H)\le n-5-r$ and $H\in
Ex(n-3;\{K_{1,n-4-r},T_{n-2-r}^*\})$.
\par Now we assume $r\ge \f{n-7}2$.
 If $\Delta(H)=n-5-r$, we may
assume $d(v_1)=n-5-r$ and $\Gamma_H(v_1)=\{v_2,\ldots,v_{n-4-r}\}$.
Since $G'$ does not contain $T_n^*$ and $d_{G'}(v_1)=n-3$, we see
that $\{v_{n-3-r},\ldots,v_{n-3}\}$ is an independent set.  As $r\le
n-6$, by the above we have $e(H)\le\sum_{i=2}^{n-4-r}d_H(v_i)\le
(n-5-r)^2.$ Since $r\ge \f{n-7}2$ we have $n-3\ge 2(n-5-r)$. Set
$H'=K_{n-5-r,n-5-r}\cup (3r+9-n)K_1$. Then $|V(H')|=n-1+r$ and
$e(H')=(n-5-r)^2$, $\Delta(H')=n-5-r$ and $H'$ does not contain
$T_{n-2-r}^*$. As $G'$ is an extremal graph, by the above we must
have $e(H)=e(H')=(n-5-r)^2$. If $\Delta(H)<n-5-r$, then clearly
$H\in Ex(n-3;K_{1,n-5-r})$. Using Theorem 2.1 we see that
$e(H)=ex(n-3;K_{1,n-5-r})=[\f{(n-3)(n-6-r)}2]$. Therefore,
$e(H)=\max\{(n-5-r)^2,[\f{(n-3)(n-6-r)}2]\}$ and so
$$\align &ex(n-1+r;\{K_{1,n-2},T_n^*\})
\\&=e(G)=e(G')=(n-3)(r+2)
+\max\Big\{(n-5-r)^2,\Big[\f{(n-3)(n-6-r)}2\Big]\Big\}.\endalign$$
This completes the proof.
\end{proof}
\pro{\noindent Theorem 4.4} Let $p,n\in\Bbb N$,
$p\ge n\ge 11$, $r\in\{2,3,\ldots,n-6\}$ and $p\e r\mod{n-1}$. Let
$m\in\{0,1,\ldots,r+1\}$ be given by $n-3\e m\mod{r+2}$. Then
$$\aligned &ex(p;T_n^*)\\&=\cases \big[\f{(n-2)(p-1)-2r-m-3}2\big]
&\t{if $r\ge 4$ and
$2\le m\le r-1$,}\\\f{(n-2)(p-1)-m(r+2-m)-r-1}2
&\t{otherwise}.\endcases\endaligned$$\endpro

%   Proof.
\noindent{\it Proof.\/}
Suppose $s=[\f{n-3}{r+2}]$. Then $n-3=s(r+2)+m$. As $r+2<n-3$ we see
that $s\in\Bbb N$. We claim that
$$\aligned &ex(n-1+r;\{K_{1,n-2},T_n^*\})
\\&=\f{(n-3-m)(n-1+r+m)}2+\max\Big\{m^2,\Big[\f{(r+2+m)(m-1)}2\Big]\Big\}.\endaligned\tag 4.2$$
When $s=1$ we have $n-5-r=m<r+2$ and so $\f{n-7}2<r<n-5$. Thus
applying Lemma 4.7 we have
$$\align
&ex(n-1+r;\{K_{1,n-2},T_n^*\})
\\&=(n-3)(r+2)+\max\Big\{(n-5-r)^2,\Big[\f{(n-6-r)(n-3)
}2\Big]\Big\}
\\&=\f{(n-3-m)(n-1+r+m)}2+\max\Big\{m^2,\Big[\f{(r+2+m)(m-1)}2\Big]\Big\}.\endalign
$$ So (4.2) holds.
\par From now on we assume $s\ge 2$.
For $i=0,1,\ldots,s-2$ we have $n-i(r+2)-5\ge n-3-(s-2)(r+2)-2\ge
2(r+2)-2>r\ge 2$. Thus, by Lemma 4.7 we have
$$\align &ex(n-3+r+2-i(r+2);\{K_{1,n-i(r+2)-2},T_{n-i(r+2)}^*\})
\\&=(r+2)(n-3-i(r+2))
 \\&\qq+ex(n-3-i(r+2);\{K_{1,n-(i+1)(r+2)-2},T_{n-(i+1)(r+2)}^*\}).\endalign$$
 Hence
 $$\align
 &ex(n-1+r;\{K_{1,n-2},T_n^*\})-ex(2(r+2)+m;\{K_{1,m+r+3},T_{m+r+5}^*\})
 \\&= ex(n-3+r+2;\{K_{1,n-2},T_n^*\})\\&\qq-ex(n-3-(s-2)(r+2);
 \{K_{1,n-(s-1)(r+2)-2},T_{n-(s-1)(r+2)}^*\})
\\&=\sum_{i=0}^{s-2}\Big(ex(n-3+r+2-i(r+2);\{K_{1,n-i(r+2)-2},T_{n-i(r+2)}^*\})
\\&\qq\qq-ex(n-3-i(r+2);\{K_{1,n-(i+1)(r+2)-2},T_{n-(i+1)(r+2)}^*\})\Big)
\\&=\sum_{i=0}^{s-2}(r+2)(n-3-i(r+2)).
\endalign$$
Set $n'=m+r+5$. As $r>m-2$ and $r\ge 2$, we have $\f{n'-7}2<r\le
n'-5$ and $n'\ge r+5\ge 7$. Thus, by Lemma 4.7 we have
$$\align &ex(2(r+2)+m;\{K_{1,m+r+3},T_{m+r+5}^*\})
\\&=ex(n'-1+r;\{K_{1,n'-2},T_{n'}^*\})
\\&=(n'-3)(r+2)+\max\Big\{(n'-5-r)^2,\Big[\f{(n'-6-r)(n'-3)}2\Big]\Big\}
\\&=(r+2)(n-3-(s-1)(r+2))+\max\Big\{m^2,\Big[\f{(m-1)(m+r+2)}2\Big]\Big\}.
\endalign$$
Therefore,
$$\align &ex(n-1+r;\{K_{1,n-2},T_n^*\})
\\&=\sum_{i=0}^{s-1}(r+2)(n-3-i(r+2))
+\max\Big\{m^2,\Big[\f{(m-1)(m+r+2)}2\Big]\Big\}.\endalign$$ As
$$\align &\sum_{i=0}^{s-1}(r+2)(n-3-i(r+2))
\\&=(r+2)\Big((n-3)s-(r+2)\f{(s-1)s}2\Big)=\f{s(r+2)}2\big(2(n-3)-(s-1)(r+2)\big)
\\&=\f{(n-3-m)(n-1+r+m)}2,\endalign$$
from the above we see that (4.2) is also true for $s\ge 2$.
\par
Observe  that $\f{(m+r+2)(m-1)}2=m^2+\f{(r-m)(m-1)-2}2$. For
$m=0,1,r,r+1,$ we have $(r-m)(m-1)-2\le 0.$ Now assume $2\le m\le
r-1$. If $r=3$, then $m=2$ and so $(r-m)(m-1)-2=-1<0$. If $r\ge 4$,
 then clearly $(r-m)(m-1)-2\ge 0.$
Thus, by (4.2) and the above we obtain
$$\aligned &ex(n-1+r;\{K_{1,n-2},T_n^*\})
\\&=\cases \f{(n-3-m)(n-1+r+m)}2+[\f{(r+2+m)(m-1)}2]\\\qq\qq\qq\qq\qq\t{if}\ r\ge 4
\ \t{and}\  2\le m\le r-1,\\\f{(n-3-m)(n-1+r+m)}2+m^2
\\\qq\qq\qq\qq\qq\t{otherwise}.\endcases\endaligned\tag 4.3$$
 \par For $r=2$ we have $m\le r+1\le 3$.
Let $G\in Ex(n+1;T_n^*).$ If $\Delta(G)\ge n-2,$
 by Lemmas 4.1 and 4.2  we have $G=K_{n-1}\cup K_2.$ Thus, $e(G)=\b{n-1}2+1.$
If $\Delta(G)\le n-3,$ then $G\in Ex(n+1;\{K_{1,n-2},T_n^*\}).$
Thus, applying (4.3) we have
$$\align &ex(n+1;T_n^*)\\&=\max\Big\{\f{(n-1)(n-2)}2+1,ex(n+1;\{K_{1,n-2},T_n^*\})\Big\}
\\&=\max\Big\{\f{(n-1)(n-2)}2+1,\f{(n-3-m)(n+1+m)}2+m^2\Big\}
\\&=\f{(n-3-m)(n+1+m)}2+m^2+\max\Big\{0,-\f{(m-2)^2+n-11}2\Big\}
\\&=\f{(n-3-m)(n+1+m)}2+m^2.
\endalign$$
  For $r\ge 3$, by Lemma 4.6 we have
$ex(n-1+r;T_n^*)=ex(n-1+r;\{K_{1,n-2},T_n^*\})$. Thus applying (4.3)
we obtain
$$\aligned &ex(n-1+r;T_n^*)
\\&=\cases \f{(n-3-m)(n-1+r+m)}2+[\f{(r+2+m)(m-1)}2]\\\qq\qq\qq\qq\qq\t{if}\ r\ge 4
\ \t{and}\  2\le m\le
r-1,\\\f{(n-3-m)(n-1+r+m)}2+m^2\\\qq\qq\qq\qq\qq\t{otherwise}.\endcases\endaligned\tag
4.4$$ By the previous argument, (4.4) is also true for $r=2$.
\par
Now suppose $p=k(n-1)+r$. Then $k\in\Bbb N$. Combining (4.4) with
Theorem 4.1 we deduce the following result:
$$\aligned &ex(p;T_n^*)\\&=\cases (k-1)\b{n-1}2+\f{(n-3-m)(n-1+r+m)}2+
\Big[\f{(r+2+m)(m-1)}2\Big]\\\qq\qq\qq\qq\qq\qq\qq\qq\qq\t{if $r\ge 4$ and
$2\le m\le r-1$,}\\(k-1)\b{n-1}2+\f{(n-3-m)(n-1+r+m)}2+m^2
\q\t{otherwise}.\endcases\endaligned$$ To see the result, we note
that
$$\align &(k-1)\b{n-1}2+\f{(n-3-m)(n-1+r+m)}2+
\Big[\f{(r+2+m)(m-1)}2\Big]
\\&=\Big[\f{(n-2)(p-1)-2r-m-3}2\Big]
\endalign$$
and
\allowdisplaybreaks\begin{gather*}\hskip+3.5cm(k-1)\b{n-1}2+\f{(n-3-m)(n-1+r+m)}2+m^2\hskip+3.5cm
\\
\hskip+3.5cm=\f{(n-2)(p-1)-m(r+2-m)-r-1}2.\hskip+3.5cm\square
\end{gather*}
\vskip+0.2cm
\pro{\noindent Corollary 4.1} Suppose $p,n,r\in\Bbb N$, $p\ge n\ge 11$,
$\f{n-7}2<r\le n-6$ and $p\e r\mod{n-1}$. Then
$$ex(p;T_n^*)=\cases \Big[\f{(n-2)(p-2)-r}2\Big]&\t{if $\f{n-4}2\le r\le
n-7$,}
\\\f{(n-2)(p-3)}2+3&\t{if $r=n-6$,}
\\\f{(n-2)(2p-5)+7}4&\t{if $r=\f{n-5}2$,}
\\\f{(n-2)(p-2)}2+1&\t{if $r=\f{n-6}2$.}
\endcases$$
\endpro
%   Proof.
\begin{proof}
Clearly $r>\f{n-7}2\ge 2$. Set $m=n-5-r$. Then $1\le m<r+2$
and $n-3\e m\mod{r+2}$. It is evident that
$$2\le m\le r-1\iff \f{n-4}2\le r\le n-7.$$
As $n\ge 11$ we see that $r\ge \f{n-4}2$ implies $r\ge 4$. Now
applying Theorem 4.4 we deduce that
$$ex(p;T_n^*)=\cases \Big[\f{(n-2)(p-1)-2r-(n-5-r)-3}2\Big]=\Big[\f{(n-2)(p-2)-r}2\Big]\\\qq\qq\qq\qq\qq\t{if $\f{n-4}2\le r\le
n-7$,}
\\\f{(n-2)(p-1)-(n-5-r)(r+2-(n-5-r))-r-1}2\\\qq\qq\qq\qq\qq\t{if
$r=n-6\  \text{or}\ [\f{n-5}2]$.}\endcases$$ This yields the result.
\end{proof}
\pro{\noindent Corollary 4.2} Suppose $p,n\in\Bbb N$, $p\ge n\ge 11$, $2\nmid
n$ and $p\e \f{n-7}2\mod{n-1}$. Then
$$ex(p;T_n^*)=\f{(n-2)(2p-3)+3}4.$$
\endpro
%   Proof.
\begin{proof}
 Taking $r=\f{n-7}2$ and $m=0$ in Theorem 4.4 we derive the
result.
\end{proof}
\pro{\noindent Corollary 4.3} Suppose $p,n\in\Bbb N$, $p\ge n\ge 11$ and
$(n-1)\mid (p-2)$. Then
$$ex(p;T_n^*)=\cases ((n-2)(p-1)-6)/2&\t{if $n\e 0\mod 2$,}
\\((n-2)(p-1)-7)/2&\t{if $n\e 1\mod 4$,}
\\((n-2)(p-1)-3)/2&\t{if $n\e 3\mod 4$.}
\endcases$$
\endpro
%   Proof.
\begin{proof}
Let $m\in\{0,1,2,3\}$ be given by $n-3\e m\mod 4$. Then
clearly $m=1,2,3\ \t{or}\ 0$ according as $n\e 0,1,2\ \t{or}\ 3\mod
4$. Now putting $r=2$ in Theorem 4.4 and applying the above we
obtain the result.
\end{proof}
\pro{\noindent Corollary 4.4} Suppose $p,n\in\Bbb N$, $p\ge n\ge 11$ and
$(n-1)\mid (p-3)$. Then
$$ex(p;T_n^*)=\cases (n-2)(p-1)/2-2&\t{if $n\e 3\mod 5$,}
\\(n-2)(p-1)/2-4&\t{if $n\e 2,4\mod 5$,}
\\(n-2)(p-1)/2-5&\t{if $n\e 0,1\mod 5$.}
\endcases$$
\endpro
%   Proof.
\begin{proof}
Let $m\in\{0,1,2,3,4\}$ be given by $n-3\e m\mod 5$. Then clearly
$m=2,3,4,0\ \t{or}\ 1$ according as $n\e 0,1,2,3\ \t{or}\ 4\mod 5$.
Now putting $r=3$ in Theorem 4.4 and applying the above we obtain
the result.\end{proof}
\par In a similar way, putting $r=4$ in Theorem 4.4 we deduce the
following result.

 \pro{\noindent Corollary 4.5} Suppose $p,n\in\Bbb N$, $p\ge
n\ge 11$ and $(n-1)\mid (p-4)$. Then
$$ex(p;T_n^*)=\cases (n-2)(p-1)/2-7&\t{if $n\e 0\mod 6$,}
\\(n-2)(p-1)/2-5&\t{if $n\e \pm 2\mod 6$,}
\\((n-2)(p-1)-13)/2&\t{if $n\e \pm 1\mod 6$,}
\\((n-2)(p-1)-5)/2&\t{if $n\e 3\mod 6$.}
\endcases$$
\endpro
 \pro{\noindent Corollary 4.6} Suppose $p\in\Bbb N$, $p\ge
11$, $r\in\{0,1,\ldots,9\}$ and $p\e r\mod{10}$. Then
$$ex(p;T_{11}^*)=\cases (9p-r(10-r))/2&\t{if $r\in\{0,1,7,8,9\}$,}
\\(9p-12)/2&\t{if $r=2$,}
\\(9p-19)/2&\t{if $r=3$,}
\\(9p-22)/2&\t{if $r=4$,}
\\(9p-21)/2&\t{if $r=5$,}
\\(9p-16)/2&\t{if $r=6$.}\endcases$$
\endpro
%   Proof.
\begin{proof}
The result follows from Theorems 4.1-4.3 and Corollaries 4.1-4.2.
\end{proof}
 \pro{\noindent Theorem 4.5}
Let $p,n\in\Bbb N$ with $6\le n\le 10$ and $p\ge n$, and let
$r\in\{0,1,\ldots,n-2\}$ be given by $p\e r\mod{n-1}$.
\par $(\t{\rm i})$ If $n=6,7$, then
$ex(p;T_n^*)=\f{(n-2)p-r(n-1-r)}2$.
\par $(\t{\rm ii})$ If $n=8,9$, then
$$ex(p;T_n^*)=\cases \f{(n-2)p-r(n-1-r)}2&\t{if $r\not=n-5$,}
\\\f{(n-2)(p-2)}2+1&\t{if $r=n-5$.}\endcases$$
\par $(\t{\rm iii})$ If $n=10$, then
$$ex(p;T_n^*)=\cases 4p-r(9-r)/2&\t{if $r\not=4,5$,}
\\4p-7&\t{if $r=5$,}\\4p-9&\t{if $r=4$.}
\endcases$$
\endpro
%   Proof.
\begin{proof}
For $r\in\{0,1,n-5,n-4,n-3,n-2\}$ the result follows from
Theorems 4.1, 4.2 and 4.3. Now assume $2\le r\le n-6$. Then $r\ge
2>\f{n-7}2$. By Lemma 4.7 we have
$$\align
&ex(n-1+r;\{K_{1,n-2},T_n^*\})\\&=(n-3)(r+2)+\max\Big\{(n-5-r)^2,\Big[\f{(n-6-r)(n-3)
}2\Big]\Big\}.\endalign$$ If $G\in Ex(n-1+r;T_n^*)$ and
$\Delta(G)\ge n-2$, using Lemmas 4.1 and 4.2 we see that $G\cong
K_{n-1}\cup K_r$. Thus,
$$\align ex(n-1+r;T_n^*)&=\max\Big\{\b{n-1}2+\b
r2,ex(n-1+r;\{K_{1,n-2},T_n^*\})\Big\}
\\&=\max\Big\{\b{n-1}2+\b
r2,(n-3)(r+2)\\&\qq+\max\Big\{(n-5-r)^2,\Big[\f{(n-6-r)(n-3)
}2\Big]\Big\}\Big\}.\endalign$$ From this we deduce that
$$\align &ex(7+2;T_8^*)=\b 72+\b 22,\q ex(8+2;T_9^*)=\b 82+\b 22,
\\&ex(8+3;T_9^*)=\b 82+\b 32,\q ex(9+2;T_{10}^*)=\b 92+\b 22,
\\&ex(9+3;T_{10}^*)=\b 92+\b 32,\q ex(9+4;T_{10}^*)=43.
\endalign$$
Suppose $p=k(n-1)+r$. Then $k\in\Bbb N$. By Theorem 4.1,
$$\align ex(p;T_n^*)&=(k-1)\b{n-1}2+ex(n-1+r;T_n^*)
\\&=\f{(n-2)(p-r)}2+ex(n-1+r;T_n^*)-\b{n-1}2.\endalign$$ Now combining all the
above we deduce the result.
\end{proof}
\section*{Acknowledgements}

The first author is supported by the Natural Sciences Foundation of China (grant No. 10971078).

\label{lastpage-01}
\end{document}